\documentstyle[12pt]{article}
\begin{document}

\title{$U_q(sl(2))$ as Dynamical Symmetry Algebra of the
Quantum Hall Effect}

\author{\"{O}mer F. Dayi\thanks{Talk presented in
       ``Quantum Groups, Deformations and  Contractions", 
       Istanbul 1997.
       E-mail: dayi@gursey.gov.tr.}\\
       Feza G\"{u}rsey Institute,\\
       P.O.Box 6, 81220 \c{C}engelk\"{o}y--Istanbul, Turkey. \\ }

\date{}
\maketitle

\begin{abstract}
Quantum Hall effect wavefunctions
corresponding to the filling factors
$\mbox{1/2p+1}$, $\mbox{2/2p+1,}\ \cdots ,\ \mbox{2p/2p+1,}\ 1,$
are shown to form a basis
of   irreducible cyclic representation
of the quantum algebra $U_q(sl(2))$ at $q^{2p+1}=1.$ 
Thus, the wavefunctions $\Psi_{P/Q}$ possessing  filling
factors $P/Q<1$ where $Q$ is odd and $P,\ Q$ are
relatively prime integers are classified in terms
of $U_q(sl(2)).$ 
Adopted as dynamical symmetry 
this leads to non--existence
of a ``universal microscopic theory"
of the quantum Hall effect,
defined as
the eigenvalue problem of
a differential operator ${\cal O}$:
${\cal O}\Psi_{\{ \nu \}}=\ell_{\{ \nu \}}\Psi_{\{ \nu \}}$
for $\nu =1,1/3,2/3,\cdots ,$
in the complex plane.
\end{abstract}

\vspace{1cm}
\noindent
{\large \bf 1. Introduction:}
\vspace{.5cm}

\noindent
Microscopic theory of the fractional
quantum Hall effect (QHE) is not well established.
Its theoretical understanding mostly is due to
trial wavefunctions \cite{qhe}.
For filling factors $1/m$ where $m$ is an odd integer, 
trial wavefunctions were given by 
Laughlin\cite{la} .
Trial wavefunctions for the other filling factors
$\nu =P/Q<1,$ where $P,\ Q$ are relatively prime integers
and $Q$ is odd,
were constructed in terms of some hierarchy
schemes\cite{owf}--\cite{jwz} 
where they were obtained from a parent
state which is a full filled Landau level or
a Laughlin wavefunction.
However, general properties of
the QHE should be independent of the explicit form of
the trial wavefunctions, but depend on their universal
features as their  orthogonality.

The integral QHE is understood in terms of non-interacting
electrons which fully fill a certain amount of Landau levels.
i.e. there exists  a microscopic
hamiltonian $h$ independent of
filling factor $n=1,2,\cdots ,$ satisfying
$h\Phi_n=E_n\Phi_n,$ where $\Phi_n$ is the wavefunction
corresponding the filling factor $n,$
and $E_n$ is the related eigenvalue.
In fact, solving  the eigenvalue problem of a given
differential operator ${\cal O}$ which defines a physical system,
is the usual procedure in
quantum mechanics. Once the eigenfunctions of ${\cal O}$
are found, they may be classified
as representations of an algebra (group)
thus named dynamical symmetry algebra of the
system. In the contrary,
if the underlying differential operator
of a physical system is not available,
knowing dynamical symmetry algebra of
the system can give some hints about it.

When one deals with the QHE $(\nu \leq 1)$
a ``universal microscopic theory" given
by a differential operator ${\cal O}$ which
is independent of the available filling factors
$\{\nu \}=1,\ 1/3,\ 2/3,\ 1/5,\cdots ,$
and satisfies the eigenvalue equation
${\cal O} \Psi_{\{\nu \}} =\ell_{\{\nu \}} \Psi_{\{\nu \}} ,$
is not known.
Microscopic theories which we know are
given for one value of the filling factor
$\nu ,$ and their excitations,
e.g. see \cite{eft} and the references therein.
i.e. they are given in terms of the eigenvalue equations
as $H_\nu \Psi_{\nu ,k} =e_{\nu ,k}\Psi_{\nu ,k},$
where $k$ labels the energy eigenvalues and the ground
state $\Psi_{\nu ,0}$  can be a full filled Landau level
or a Laughlin wavefunction.
These are effective theories which depend on
the number of the  levels  occupied
by the electrons.

We utilize orthogonality of the QHE states for
different filling factors, independent of their
explicit form, to
show that they
can be classified
as irreducible cyclic representations
of $U_q(sl(2))$ at  roots of unity\cite{om}.
In our scheme,  states 
corresponding to 
filling factors  possessing a common denominator 
are in the same representation.
Based on this classification
 $U_q(sl(2))$ at roots of unity is proposed to be
dynamical symmetry algebra of the QHE, which
interrelates states of different filling
factors. 
This leads
to the conclusion that a ``universal microscopic theory"
of the fractional QHE, 
in the common sense, does not exist. 
This would be the explanation why such a microscopic theory
of the fractional QHE is not known.

\vspace{1cm}
\noindent
{\large \bf 2. Cyclic Representation of $U_q(sl(2)):$}
\vspace{.5cm}
\noindent

The deformed algebra $U_q(sl(2))$  
\begin{eqnarray}
[E_+,E_-] & = & \frac{K-K^{-1}}{q-q^{-1}},  \nonumber \\
KE_{\pm}K^{-1} & = & q^{\pm 2}E_{\pm}. \label{alg}
\end{eqnarray}
at roots of unity i.e. $q^{2p+1}=1,$ $p$ a positive integer, 
has a finite dimensional irreducible representation 
which has no classical finite dimensional analog.
This is the cyclic representation whose dimension 
is $2p+1$\cite{rep}. Cyclic means that there 
are no heighest or lowest weight states in the spectrum.
i.e. $E_+|\cdots >\neq 0$ and $E_-|\cdots >\neq 0$
for any state.

When $q^{2p+1}=1$ irreducible cyclic representation
of $U_q(sl(2))$
can be written in some basis
$\{ v_0,v_1, \cdots, v_{2p} \}$ as
\begin{eqnarray}
Kv_m & = & \lambda q^{-2m} v_m, \nonumber \\
E_+v_m & = & g_m v_{m+1} , \label{rep} \\
E_-v_m & = & f_m v_{m-1} , \nonumber
\end{eqnarray}
where $m=0,\cdots ,2p,$ and we defined $v_0\equiv v_{2p+1},\
v_{-1}\equiv v_{2p}.$ $\lambda,$ $g_m$ , and
$f_m$ are some complex constants which are nonzero and
in the case of requesting 
that the representation in
unitary, we should restrict their values
such that
\begin{equation}
\label{uni}
K^\dagger =K^{-1};\   E^\dagger_-=E_+.
\end{equation}

\vspace{1cm}
\noindent
{\large \bf 3. Classification:}
\vspace{.5cm}

\noindent
QHE trial wavefunctions
in the standard hierarchy scheme
are given by\cite{owf},\cite{re}
\begin{eqnarray}
\psi_\nu (z_1,\cdots ,z_{N_0})&  = & \int
\prod_{\alpha =1}^{r}
\prod^{N_\alpha }_{i_\alpha =1}[d^2z^{(\alpha )}_{i_\alpha}]
e^{-\frac{1}{2}\sum_1^{N_0} |z_k|^2}
\prod_{\beta =0}^{r}
\prod_{i_\beta <j_\beta}^{N_\beta} 
(z^{(\beta )}_{i_\beta}-z^{(\beta )}_{j_\beta})^{a_\beta} \nonumber \\
& & \times \prod_{i_{\beta+1},j_\beta =1}^{N_{\beta+1},N_\beta}
(z^{(\beta +1 )}_{i_{\beta +1}}-
z^{(\beta )}_{j_\beta})^{b_{\beta ,\beta +1}},   \label{sth}
\end{eqnarray}
where 
$z_{i_0}^{(0)} \equiv z_i .$
The measure 
$\prod [d^2z^{(\alpha )}_{i_\alpha}] $
depends on $a_\beta$ and 
$|z^{(\beta )}_{i_\beta}-z^{(\beta )}_{j_\beta}|,$ however the 
detailed form of it 
does not affect the filling factor $\nu =P/Q .$
$a_0$ is an odd positive integer,
$a_\alpha $ for $\alpha \neq 0$ are even integers
which can be positive or negative and
$b_{\beta +1, \beta}=\pm 1,$ except 
$b_{r,r+1}=0.$
By placing the $N_0$ electrons on a spherical
surface in a monopole magnetic field, one can
find that filling factor of (\ref{sth}) is given by
\begin{equation}
\label{gff}
\nu =\frac{1}{a_0 - \frac{1}{a_1- \frac{1}{\cdots -\frac{1}{a_r}}}}.
\end{equation}
Factors with negative powers may be replaced 
by complex--conjugate factors with positive powers
multiplied by some exponential factors.
Hence, (\ref{sth}) can equivalently be given as\cite{bw}
\begin{eqnarray}
\psi_\nu (z_1,\cdots ,z_{N_0})& = & \int
\prod_{\alpha =1}^{r}
\left[ \prod_{i_\alpha =1}^{N_\alpha} d^2z^{(\alpha )}_{i_\alpha}
\prod_{i_\alpha <j_\alpha}^{N_\alpha} 
|z^{(\alpha )}_{i_\alpha}
-z^{(\alpha )}_{j_\alpha}|^{2(-1)^{\alpha}\theta_\alpha} 
e^{-|q_\alpha |\sum_{i_\alpha}
|z^{(\alpha )}_{i_\alpha}|^2} \right] \nonumber \\
& & \times e^{-\frac{1}{2}\sum_1^{N_0} |z_{k}|^2}
\prod_{\beta =0}^{r} \prod_{i_\beta < j_\beta}^{N_\beta}
(\tilde{z}^{(\beta )}_{i_\beta} 
-\tilde{z}^{(\beta )}_{j_\beta})^{p_\beta}
\prod_{i_{\beta+1},j_\beta=1}^{N_{\beta+1},N_\beta}
(\bar{\tilde{z}}_{i_{\beta +1}}^{(\beta +1) }
-\tilde{z}^{(\beta )}_{j_\beta}), \label{bws}
\end{eqnarray}
where $\tilde{z}^{(\beta)}_{i_\beta} =z^{(\beta )}_{i_\beta}$ 
for $\beta =$ even
and $\tilde{z}^{(\beta)}_{i_\beta} =\bar{z}^{(\beta)}_{i_\beta}$ 
for $\beta =$ odd and
\begin{eqnarray}
\theta_0=0,  & \theta_r  = 
\frac{(-1)^r}{p_{r-1}-(-1)^r\theta_{r-1}}, \nonumber \\
q_0  =-1, & q_r  =   (-1)^{r+1}q_{r-1}\theta_r. \nonumber
\end{eqnarray}
Now, filling factor is
\begin{equation}
\label{ffbw}
\nu =\frac{1}{p_0 +\frac{1}{p_1+ \frac{1}{\cdots +\frac{1}{p_r}}}},
\end{equation}
where $p_0$ is odd and the other $p_i$ are even integers.

By generalizing the calculations of 
Laughlin given in Ref. \cite{qhe}
and making use of the scalar product
\begin{equation}
\label{spr}
(\psi_\nu , \psi_{\nu^\prime})\equiv
\int d^2z_1\cdots d^2z_N \bar{\psi}_\nu (z_1\cdots ,z_N) 
\psi_{\nu^\prime} (z_1\cdots ,z_N) .
\end{equation}
one can show that $\psi_\nu$ states are orthogonal\cite{re}.

To emphasize the second quantized character of our construction
let us introduce the states
\begin{equation}
\label{gb}
|i,\ p>_T  =  \int d^2z_{1}\cdots d^2z_{N_0}
e^{- \frac{1}{2}
\sum_{k=1}^{N_0}|z_k|^2}
\psi_{\frac{i}{2p+1}}(z_1,\cdots ,z_{N_0})
|z_1,\cdots ,z_{N_0}>,
\end{equation}
where $i=1,\cdots ,2p+1;\  p=1,2, \cdots ,$  
so that any filling factor $\nu =P/Q$ is considered.
We used the vectors
\begin{equation}
\label{ve}
|z_1,\cdots ,z_{N_0}> = \frac{1}{\sqrt{N_0!}}
\varphi^\dagger (z_1) \cdots 
\varphi^\dagger (z_{N_0}) |0>.
\end{equation}
The fermionic operators $\varphi (z),\ \varphi^\dagger (z)$ 
satisfy the anticommutation relation
\[
\{\varphi^\dagger (z), \varphi (z^\prime )\} =e^{z^\prime \bar{z}}.
\]
The subscript $T$ denotes the fact that trial wave
functions are used to give an explicit realization.
The states (\ref{gb}) are orthonormal:
\begin{equation}
\label{ort}
_T<i,\ p|j,\ p^\prime >_T=\delta_{i,j}\delta_{p,p^\prime}.
\end{equation}

We have shown that the states
$|i,\ p>_T $ are orthonormal by using the explicit form
of trial wavefunctions. However,
this should be a universal feature of QHE wavefunctions.
Then,
even if we do not know the explicit form, we can say that
exact states of the QHE  which we  indicate with $|i,\ p>,$
should be orthonormal:
\begin{equation}
\label{ort1}
<i,\ p|j,\ p^\prime >=\delta_{i,j}\delta_{p,p^\prime}.
\end{equation}
Indeed, in the following we will use this universal
property of QHE states without referring to any
trial wavefunction.

Let us deal with the states
\begin{equation}
\label{hik}
|1,\ p>,\ |2,\ p>,\cdots,\ |2p,\ p>,\ |2p+1,\ p>,
\end{equation}
corresponding to the filling factors
\begin{equation}
\label{gnu}
\nu =\frac{1}{2p+1}, \frac{2}{2p+1},\cdots ,\frac{2p}{2p+1},1.
\end{equation}

Define the following
second quantized operators acting in the space
spanned by the states  (\ref{hik}),
\begin{eqnarray}
\tilde{K}
& = & \sum_{i=1}^{2p+1} q^i|i,\ p>
<i,\ p| ,\label {gc0} \\
\tilde{E}_+
& = & \sum_{i=1}^{2p+1} a_i|i,\ p>
<i+2,\ p| , \label{gc1} \\
\tilde{E}_-
&= & \sum_{i=1}^{2p+1} \bar{a}_i
|i+2,\ p><i,\ p| , \label{gc2}
\end{eqnarray}
where 
\begin{equation}
\label{root}
q^{2p+1} =1.
\end{equation}
To obtain the compact forms we adopted the definitions
\[
|2p+2,\ p>\equiv |1,\ p> ,\  |2p+3,\ p>\equiv |2,\ p>.
\]

By using the orthonormality condition (\ref{ort})
one observes that
inverse of $\tilde{K}$ is 
\begin{equation}
\label{k-1}
{\tilde{K}}^{-1} =  \sum_{i=1}^{2p+1} q^{-i}|i,\ p><i,\ p| =
{\tilde{K}}^\dagger .
\end{equation}

Let the coefficients $a_i$ 
are nonzero and satisfy
\begin{eqnarray*}
|a_{2p+1}|^2  - |a_{2p-1}|^2 & = & 0, \\
|a_{2p}|^2 -  |a_{2p-2}|^2 & = & -1, \\
|a_{l+2}|^2  -  |a_{l}|^2 & = &
\frac{q^{l+2}-q^{-l-2}}{q-q^{-1}}, \\
\end{eqnarray*}
where $l=-1,0,\cdots (2p-3);\ a_{-1}\equiv a_{2p},\
a_0\equiv a_{2p+1}.$ 
Then,
in terms of the basis $(|1,\ p>,\ \cdots ,\ |2p+1,\ p>)$
the operators
(\ref{gc0})--(\ref{gc2}) lead to a $(2p+1)$
dimensional unitary irreducible cyclic  representation
of $U_q(sl(2))$ at q satisfying (\ref{root}).

\vspace{1cm}
\noindent
{\large  \bf 4. Discussions:}
\vspace{.5cm}

\noindent
It is shown that QHE wavefunctions can be classified as
irreducible cyclic representations of
$U_q(sl(2))$ at roots of unity in a very natural way.
This naturalness follows from the fact that
the most significant physical quantity
of the QHE $\nu =P/Q $ 
fits very well with the integer ($m$ in (\ref{rep}))
characterizing irreducible cyclic representations of 
$U_q(sl(2)).$  

Because of this classification we can propose 
$U_q(sl(2))$ at  roots of unity
as dynamical symmetry algebra of the QHE.
Obviously,  any set of orthogonal
states possessing a quantum number
which permits a partition of 
unity like $\nu ,$ 
\[
\sum_{i=1}^{2p+1}\frac{ \nu (|i,\ 2p+1>)}{p+1} =1 ,
\]
can be classified 
as irreducible cyclic  representation
of $U_q(sl(2))$ at a root of unity. 
Hence,
our observation is not enough to prove that
$U_q(sl(2))$ at roots of unity
is the real dynamical symmetry algebra of the QHE.
However, there exists another evidence to believe that
the proposed dynamical symmetry is the one chosen
by nature. If there exists a ``universal 
microscopic theory" of the QHE given in terms
of a differential operator depending
on $z_k,\ \bar{z}_k$ and their derivatives
and moreover, possessing this dynamical symmetry,
it should be in the form
$$
{\cal O}_{\tilde{q}} \Phi_i =\ell_i \Phi_i;\ i=1,\cdots ,2p+1;\
{\tilde{q}}^{2p+1}=1.
$$
Here, ${\cal O}_{\tilde{q}}$ denotes a differential operator
which is a function of the generators of the dynamical
symmetry algebra $U_q(sl(2))$ and
$\Phi_i$ are its eigenfunctions corresponding to the
eigenvalues $\ell_i$
and possessing the filling factors
$\nu (\Phi_i)=i/2p+1.$
Thus,
one should find a differential realization of
the generators of
$U_q(sl(2))$ which
determine the cyclic representations.
These differential operators should act on 
polynomials in $z_k.$ Because, the
$\nu =1$ wavefunction which is an element
of the basis of
irreducible cyclic representation is exact and it is a
polynomial in $z_k.$
There are some differential realizations of these
generators leading to the 
cyclic representations of $U_q(sl(2))$ given 
in a space of polynomials if
the following equivalence relation is satisfied\cite{rep}
$$
z^{j+2p+1} \sim z^j  .
$$
But, these constraints are not permitted in the
complex plane where QHE wavefunctions should be constructed.
This leads to the conclusion that a 
``universal microscopic theory"  for the fractional QHE
(in the common sense),
does not exist if its dynamical algebra is 
$U_q(sl(2))$ at  roots of unity. This may explain
why one could not find a universal microscopic theory 
of the fractional QHE. 

How one can utilize the proposed dynamical symmetry for
the QHE to calculate some physical quantities?
Here, one of the most significant physical quantities is
the partition function which may be obtained if the
Green function on the space defined by
$U_q(sl(2))$ at q roots of unity with cyclic representation
is available. In Ref. \cite{ad} Green function on the space
defined by the q--deformed group $SU_q(2)/U(1)$
for $q$ not a root of unity is obtained
without referring to explicit form of the representations
but depending only on their general features.
We hope that a similar calculation can be used
in our case. Then, we can obtain Green function and the
partition function. This may  lead to a
decisive answer if the proposed 
dynamical symmetry is the real symmetry of the QHE
and moreover, it may give some hints about its
physical interpretation.

\end{document}